\input amstex
\documentstyle{amsppt}
\nologo
\NoBlackBoxes
\magnification = 1200
\pageheight {19cm}
\footnote {The author was partially supported by European Science Project AGE
'Algebraic Geometry in Europe',  contractn. ERBCHRXT940557 and by the Italian
Research Program  'Geometria Algebrica, Algebra Commutativa e aspetti
computazionali'$\quad \quad \quad \quad \quad \quad \quad$} \par \noindent
\topmatter
\title THE DEGREE OF THE GAUSS MAP FOR A GENERAL PRYM THETA-DIVISOR
\endtitle
\author {Alessandro Verra} 
\endauthor
\affil {Dipartimento di Matematica , Universit\'a di Roma 3 }
\endaffil
\address {Largo S. Leonardo Murialdo 1 , 00144 ROMA , Italy} 
\endaddress
\endtopmatter
\document
\bigskip \noindent
\bf 1.Introduction \rm \par \noindent An interesting approach to the theta
divisor $\Theta$ of a principally polarized abelian variety $A$ is certainly offered
by its Gauss map $\gamma: \Theta
\to \bold P^{d-1}$, $d = dim A$. This is the natural map which associates to
$x \in \Theta-Sing \Theta$ the tangent space $T_{\Theta,x} \subset H = T_{A,x}$.
Since $H$ does not depend on $x$, $T_{\Theta,x}$ is a point
of~$\bold  P(H^*) = \bold P^{d-1}$. As is well known $\gamma$ has finite
degree unless $A$ is a product, moreover $deg \gamma = d!$ if $\Theta$ is smooth.
The degree of $\gamma$ varies according to the singularities of $\Theta$ and so it
is interestingly related to the Andreotti-Mayer stratification of the moduli space of
$(A,\Theta)$. On the other hand, if we consider pairs $(A,\Theta)$ such that
$dim Sing \Theta > 0$, we realize that, essentially, the map $\gamma$ is
completely understood in the only case of Jacobians. Here
$\gamma$ has degree $\binom {2d-2}{d-1}$ and its description is encoded in the
beautiful geometry of the canonical curve. \par \noindent 
After Jacobians of curves, the next usual step consists in considering Prym varieties.
In this note we afford part of this step: we compute the degree of the Gauss map for
the theta divisor  $\Xi$ of a general Prym variety $P$. \par \noindent Assume  $P$ is
defined by a pair $(C,\eta)$, where $C$ is a smooth, irreducible, canonically
embedded curve of genus $g = d+1$ and $\eta$ is a non trivial $2$-torsion element of
$Pic^0(C)$. Then $H = H^0(\omega_C \otimes
\eta)^*$ so that $\bold P^{d-1} = \mid \omega_C \otimes \eta \mid$. Our result relies
on a geometric description we give for the map~$\gamma: \Xi \to
\mid \omega_C \otimes \eta \mid$. This relates $\Xi$ to the family of rank three
quadrics touching $C$ along a Prym-canonical divisor, that is to quadrics $q$ such
that
$q \cdot C = 2d$ with~$d \in \mid \omega_C \otimes \eta \mid$. As we will see, there
is a rational map $\alpha: \Xi \to Y^+$, where $Y^+$ is an irreducible component of
this family of quadrics. $\alpha$ is the quotient map $\Xi
\to
\Xi/<-1>$ followed by a birational isomorphism. Moreover there is also the natural
map 
$$
\lambda^+: Y^+ \to \mid \omega_C \otimes \eta \mid
$$ such that $\lambda(q) = d$. We will prove that $\gamma = \lambda \cdot \alpha$.
Then we will use the construction to show the main result of the paper: \bigskip
\noindent
\proclaim {THEOREM } Let $g \geq 3$ and let $\Xi$ be the Theta-Divisor of a
general Prym variety of dimension $g-1$. Then the degree of the Gauss map $\gamma:
\Xi \to \bold P^{g-2}$ is
$$ D(g) + 2^{g-3},
$$ where $D(g)$ is the degree of the variety of all quadrics of rank $\leq 3$ in
$\bold P^{g-1}$. \endproclaim  
\noindent  The number $D(g)$ can be found in the
classical book of Room ([Ro] p.133, cfr.~[HT]),
$$ D(g) = \frac {\binom g{g-3} \binom {g+1}{g-4} \cdot \cdot \cdot \binom {2g-5}2
\binom {2g-4}1} {\binom {2g-7}{g-3}\binom{2g-9}{g-4} \cdot \cdot \cdot \binom 32
\binom 11}.
$$ Since a general principally polarized abelian variety of dimension
$\leq 5$ is a Prym, (cfr. [B1]), our formula yelds $(g-1)! = D(g) + 2^{g-3}$ if $3
\leq g \leq 6$. If $dim P = 6$ we know from Debarre that the singular locus of
$\Xi$ consists of $16$ ordinary double points, ([D]). This implies $deg \gamma =
(g-1)!-32$ if $g=7$ and gives $(g-1)! - 32 = D(g) + 2^{g-3} = 688$.
Assuming $P$ general, these seem to be the only cases where the degree of
$\gamma$ was previously known. \par \noindent
Now we describe our strategy to compute the degree of the map 
$\lambda^+: Y^+ \to
\mid \omega_C \otimes \eta \mid$, hence to show the theorem. Assume $C$ has
general moduli, we explain in section $2$ that $Y^+$ is only one irreducible
component of the variety $Y$ of all  quadrics $q$ of  rank $\leq 3$ satisfying
$q \cdot C = 2d$, with $d \in \mid \omega_C \otimes \eta \mid$. There is another
irreducible component $Y^-$ such that $Y = Y^+ \cup Y^-$. We can consider as
above the map
$$
\lambda: Y \to \mid \omega_C \otimes \eta \mid
$$ sending $q$ to $\lambda(q) = d$, in particular $\lambda/Y^+ = \lambda^+$. 
Let $\bold Q = \bold PH^0(\Cal O_{\bold P^{g-1}}(2))$, it turns out that
$\lambda = l/Y$ where
$$ l: \bold Q \to \mid \omega_C^{\otimes 2} \mid
$$
is the restriction map. Note that $l$ is a linear projection, its center is the linear
system $I_C$ of all quadrics through $C$, which is canonically embedded in
$\bold P^{g-1}$. Let $$ \bold
Q^3 =
\lbrace q
\in
\bold Q/ rk(q) \leq 3 \rbrace,$$ restricting $l$ to $\bold Q^3$ we obtain a finite
covering of degree
$D(g) = deg \bold Q^3$. This happens because $dim \bold Q^3 = dim \mid
\omega_C^{\otimes 2} \mid$ and $I_C \cap \bold Q^3$ is empty. On the  other hand it
holds
$$ Y = l^{-1}(V) \cdot \bold Q^3,
$$ where $V$ is the embedding of $\mid \omega_C \otimes \eta \mid$ in
$\mid \omega_C^{\otimes 2} \mid$ under the map $d \to 2d$. This implies
$$ 
deg \lambda^+ + deg \lambda^- = D(g),
$$ 
where $\lambda^- = \lambda/Y^-$. To distinguish $deg \lambda^+$
from $deg \lambda^-$ we specialize $C$ to a trigonal curve. Here some nice geometry
appears, relating $Y^+$ and $Y^-$ (which are now reducible) to  theta-characteristics
on hyperelliptic curves of genus
$g-3$. \par \noindent
If $C$ is trigonal the degree of $\lambda$ is no longer $D(g)$ because
$S = I_C \cap \bold Q^3$ is not empty. Let $p = g-3$, we prove that $deg \lambda =
2^{2p}$ in this case. Hence the contribution of the excess intersection $S$ to the
degree of a general
$\lambda$ is $D(g)-2^{2p}$. \par \noindent
In addition we show that $deg \lambda^+ = 2^{p-1}(2^p+1)$ and
$deg \lambda^- = 2^{p-1}(2^p-1)$, which are the numbers of even
and odd theta's on a curve of genus $p$. The reason for these numbers is that a
trigonal $C$ is contained in a rational normal scroll $R$. Given $d \in \mid \omega_C
\otimes \eta \mid$, there exists exactly one curve $E
\in \mid \Cal O_R(2) \mid$ such that $E \cdot C = 2d$. The fibre at $d$ of
$\lambda^+$ (of $\lambda^-$) is then naturally bijective to the set of even (odd)
theta's on $E$. \par \noindent
Finally,  using a standard parity lemma in the theory of Prym varieties and its
geometric consequences, we deduce that the contribution of $S$ to the degrees of
general
$\lambda^+$ and $\lambda^-$ is the same, (cfr. 3.5(2) and lemma 3.20). Then we
can compute
$deg \lambda^+ = \frac 12 D(g)+2^{p-1}$ and $deg \lambda^- = \frac 12D(g)-2^{p-1}$.
Hence the degree of $\gamma$ is $2deg \lambda^+ = D(g)+2^{g-3}$. \par \noindent 
The idea of using rank three quadrics is in some sense related to Tjurin's paper [T],
where $\gamma$ has a complementary description. Let $\pi:
\tilde C \to C$ be the
\'etale double covering defined by $\eta$. Then $P
\subset Pic^{2g-2}(\tilde C)$ and each point $x \in \Xi - Sing \Xi$ is a
quadratic singularity for the theta divisor $\tilde \Theta = \lbrace
\tilde L \in Pic^{2g-2}(\tilde C) / h^0(\tilde L) \geq 1 \rbrace$. As is well
known the projectivized tangent cone of $\tilde \Theta$ at $x$ is a
quadric of rank $\leq 4$
$$
\tilde Q \subset \bold P(H^0(\omega_C\otimes \eta)^*\oplus H^0(\omega_C)^*).
$$ As in $[T]$ we can restrict $\tilde Q$ to $P(H^0(\omega_C\otimes
\eta)^*$: this is the hyperplane $\bold PT_{\Xi,x}$ counted twice, hence the point
$\gamma(x)$. As in section $2$ we can restrict $\tilde Q$ to $\bold PH^0(\omega_C)^*$:
this is a quadric $q$ in the family $Y^+$. \par
\noindent To avoid any problem we work on the complex field. Nevertheless it
seems quite possible that the proof extends to any algebraically closed field $k$,
$char k \neq 2$.  Finally, we wish to thank some friends for various interesting
discussions on the subject,  in particular V. Kanev, E. Sernesi, C. Ciliberto. 
\bigskip \noindent \it Some frequently used notations: \rm - $C^{(d)}$: the
$d$-symmetric product of a curve
$C$. \par \noindent
- $x+{C^{(d-1)}}$: the ample divisor $\lbrace d \in C^{(d)}/ d=x+d'\rbrace$.
\par \noindent
- $\mid L \mid$ ($\mid V \mid$): the linear system of the
line bundle $L$ (of the space $V \subset H^0(L)$). \par
\noindent - $<S>$: the linear span of a set $S$ in a vector or projective space.
\par \noindent - $X \cdot Y$: the intersection scheme
of $X$ and $Y$. $XY$: the intersection number. 
\par \noindent - $[Z]$ the class of $Z \subset X$ in the numerical equivalence ring
 of $X$. \par \noindent
 - $f^*(o)$: the scheme theoretic fibre at $o$ of a morphism $f: X \to Y$.
 \par \noindent
 - $\Cal S_t$: the fibre of a $T$-scheme $\Cal S$ at a point $t$ of $T$.
 \par \noindent
 -$E^*$: the dual of the vector (projective) space $E$. $\bold PE$: $Proj E^*$.
\bigskip \noindent \bf 2. Gauss map and rank three quadrics. \rm  \par \noindent  In
the following we introduce the basic construction relating the Gauss map to rank three
quadrics, moreover we fix our notations. The main object to be considered is
an
\'etale double covering 
$$
\pi: \tilde C \to C, 
$$ defined by the non trivial line bundle $\eta$. The corresponding involution on
$\tilde C$ is denoted by $i$. For simplicity, it will be assumed that both $\tilde C$
and $C$ are non hyperelliptic and  that $\omega_C \otimes \eta$ is very ample. The
latter property is satisfied unless $\eta \cong$  $\Cal O_C(x_1+x_2-y_1-y_2)$ with
$x_1,x_2,y_1,y_2 \in C$, ([CD], 0.6). An element $d \in \mid \omega_C \otimes \eta
\mid$ is said to be a \it Prym-canonical divisor. \rm As usual
$\pi$ induces the Norm map
$$  Nm: \tilde J \to J, 
$$  
where $\tilde J = Pic^{2g-2}(\tilde C)$ and  $J = Pic^{2g-2}(C)$. By definition the
image of $\Cal O_{\tilde C}(\Sigma x_i)$ by $Nm$ is $\Cal O_C(\Sigma \pi(x_i))$. The
fibre of $Nm$ is split in two isomorphic connected components.  In particular
$Nm^{-1}(\omega_C) = P \cup P^-$, where
$$ \text { $P = \lbrace \tilde L \in Nm^{-1}(\omega_C) / h^0(\tilde L)$ is even
$\rbrace$, $\quad P^- =
\lbrace
\tilde L \in Nm^{-1}(\omega_C) / h^0(\tilde L)$ is odd $\rbrace$}. 
$$
As usual, we will say that $P$ is \it the Prym variety of $\pi$. \rm $P$ is
biregular to a ($g-1$)-dimensional abelian variety, namely to $Im (1-i^*) \subset
Pic^0(\tilde C)$.
\rm 
The restriction to $P$ of the theta-divisor $$ \tilde \Theta = \lbrace
\tilde L\in\tilde J/ h^0(\tilde L) >0 \rbrace $$ 
is twice a principal polarization $\Xi \subset P$, in particular it holds $\Xi = Sing
\tilde \Theta \cdot P$, (cfr. [M] and [ACGH] p.295). By definition $\Xi$ is \it the
Prym Theta-divisor of $\pi$. \rm To introduce some projective geometry related to
$\Xi$ we consider the vector spaces 
$$  H = H^0(\omega_{\tilde C})^*, \quad H^-, \quad H^+ 
$$  where $H^-$ , $H^+$ are the eigenspaces of the involution induced by $i$ on $H$.
We notice some canonical identifications induced by $\pi^*$:
$H^+ = H^0(\omega_C)^*$, $H^- = H^0(\omega_C \otimes \eta)^*$. For the associated
projective spaces the notations will be
$$
\bold P = \bold PH \quad , \quad \bold P^+ = \bold PH^+ \quad , \quad \bold P^- =
\bold PH^-. \tag 2.1
$$  Let 
$$  h^+: \bold P \to \bold P^+ \quad \text {and} \quad h^-: \bold P \to \bold P^- 
$$ be the linear projections of centers $\bold P^-$ and $\bold P^+$. We have the
commutative diagram
$$
\CD {\bold P^-} @<{h^-}<< {\bold P} @>{h^+}>> {\bold P^+} \\ {\cup} @. {\cup}
@. {\cup} \\ {C^-} @<{h^-/\tilde C}<< {\tilde C} @>{h^+/ \tilde C}>> {C,} \\
\endCD \tag 2.2
$$ where $\tilde C$ is canonically embedded in $\bold P$ and 
$$  C^+ = h^+(\tilde C) \quad , \quad C^- = h^-(\tilde C). \tag 2.3
$$  It turns out that $\pi = h^+/\tilde C = h^-/\tilde C$. Moreover $\Cal O_{C^-}(1)
\cong \omega_C \otimes \eta$ so that $C^-$ is the Prym canonical embedding of $C$. On
the other hand $C^+$ is the canonical model of $C$. For simplicity we will put
$$  C = C^+, $$  
when no confusion arises. To simplify notations we put  
$$
\bold Q = \bold PH^0(\Cal O_{\bold P^+}(2)) \quad , \quad \bold B = \bold
PH^0(\omega_C^{\otimes 2}). \tag 2.4
$$
The natural restriction map will be denoted by
$$
\lambda: \bold Q \to \bold B. \tag 2.5
$$ 
$\lambda$ is induced from the standard exact sequence
$$ 0 \to H^0(\Cal I(2)) \to H^0(\Cal O_{\bold P^+}(2)) \to H^0(\omega_C^{\otimes 2})
\to 0,
$$ where $\Cal I$ is the Ideal of $C$, so $\lambda$ is the linear projection
of center $\mid \Cal I_C(2)\mid$. Let 
$$
\bold Q^r = \lbrace q \in \bold Q / rank (q) \leq r \rbrace,
$$
we will be specially interested to the variety $\bold Q^3$. \it The  degree of
$\bold Q^3$ is denoted as $D(g)$. \rm $\bold Q^3$ has dimension $3g-4$, (cfr.
[ACGH], p.100), as well as $\bold B$. Therefore
$$
\lambda/\bold Q^3: \bold Q^3 \to \bold B \tag 2.6
$$ 
is a finite morphism of degree $D(g)$ if $\mid \Cal I_C(2) \mid \cap \bold Q^3 =
\emptyset$. \it This condition is satisfied if $C$ is sufficiently general, \rm (see
remark 2.9-1). Finally let 
$$
v: \mid \omega_C \otimes \eta \mid \to \bold B. 
$$
be the 'squaring map' sending $d \in \mid \omega_C \otimes \eta \mid$ to $2d$. $v$ is
obtained from the diagram 
$$
\CD {H^0(\omega_C \otimes \eta)} @>{\sigma}>> {Sym^2H^0(\omega_C \otimes \eta)}
@>{\mu}>> {H^0(\omega_C^{\otimes 2}})\\
\endCD
$$ as the projectivization of $\mu \cdot \sigma$, where $\mu$ is the multiplication
map and $\sigma(x) = x^2$.  We leave as an exercise to show that $v$ is an embedding.
If $\mu$ is surjective then $v$ is the $2$-Veronese embedding of $\bold P^{g-2}$ in
a subspace of $\bold B$ of appropriate dimension. It is known that $\mu$ is surjective
if $C$ is general of genus $g \geq 6$, ([B2]). The image of $v$ in $\bold B$ will be
denoted as 
$$ 
V. \tag 2.7
$$
The inverse image of $V$ by $\lambda$ will be
$$ 
W. \tag 2.8
$$ 
$W$ is a cone over $V$ in the projective space $\bold Q$, in particular $Sing W = \mid
\Cal I_C(2)\mid$. \bigskip \noindent 
\bf 2.9 REMARKS \rm (1) A rank three quadric $q$ containing $C$ defines a line
bundle $L$ of  degree $\leq g-1$ with non injective Petri map. Hence $C$ is not
general.  To see this consider the moving part $M$ of the pencil which is cut on $C$
by the ruling of maximal linear subspaces of $q$. Observe that $M$ is defined by a
vector space $U \subseteq H^0(L)$, where $L$ is a line bundle satisfying  $L^{\otimes
2} \cong \omega_C(-b)$ and $b\geq 0$. By the base-point-free-pencil-trick, ([ACGH]
p.126), the Petri map $\mu: H^0(L) \otimes H^0(L(b)) \to H^0(\omega_C)$ is not
injective on $U
\otimes H^0(L(b))$. \par
\noindent  (2) For each quadric $q \in W-I_C$, it holds $q \cdot C = 2d$, where
$\lambda(q) = 2d$. Let $I_{2d}$ be the vector space of quadratic forms vanishing on
$2d$. Then $\mid I_{2d} \mid$ is a space of maximal dimension in the cone $W$. For a
general $2d$ the fibre of $\lambda/\bold Q^3$ over $2d$ is the locus of rank three
quadrics in $\mid I_{2d} \mid-I_C$. \bigskip \noindent
\proclaim {2.10 DEFINITION} The intersection scheme $Y = W \cdot
\bold Q^3$ is the variety of Prym-canonical rank three quadrics. \endproclaim
\proclaim {2.11 LEMMA} (i) Every irreducible component of $Y$ has dimension $\geq
g-2$. \par \noindent (ii) If $C$ is sufficently general $Y$ is reduced of pure
dimension $g-2$.
\endproclaim
\demo {Proof} (i) Dimension count, since $codim_{\bold Q}W = 2g-2$. (ii) If $C$ is
general the map $\lambda/\bold Q^3$ is finite, hence $Y = \lambda^*(V)$ is pure of
dimension $g-2$. In particular each irreducible component $Z$ of $Y$ is horizontal
i.e. $\lambda(Z)=V$. One can show that, moving $W$ in a smooth family, the limit of a
horizontal irreducible component is still horizontal. Then, to prove the reducedness
of a general $Y$, it suffices to produce \it one \rm $Y$ such that all its
irreducible horizontal components are reduced. This is the case if $C$ is a general
trigonal curve: see corollary 4.8.
\enddemo \noindent
Throughout the paper the notions of \it general curve \rm or of \it general cover
$\pi: \tilde C \to C$ \rm are used in different ways. The next definition will 
help to distinguish among them.
\proclaim {2.12 DEFINITION}(1) $\pi$ is Petri general if no element of $Y$
contains $C$. \par \noindent 
(2) $\pi$ is standard if $X$ and $Y$ are reduced. \par \noindent
$R_{Petri}$ and $R_{standard}$ denote the corresponding subsets in the moduli space
of $\pi$.
\endproclaim \noindent
We will see very soon that $Y$ is strongly related to the Prym-Theta divisor~$\Xi$ and
to its Gauss map.  To study this relation let
$$
\pi^{(2g-2)}: \tilde C^{(2g-2)} \to C^{(2g-2)}
$$ 
be the morphism defined by pushing-forward divisors via $\pi$ and let
$$
\Cal G \quad \text {and} \quad \Cal B
$$
be the following objects: \par \noindent
- $\Cal B$ is the Hilbert scheme of conics in $C^{(2g-2)}$, that is of connected
curves $B$ such that $B \cdot (x+C^{(2g-3)}) = 2$ and $p_a(B)=0$. 
\par \noindent
- $\Cal G = G^1_{2g-2}(\tilde C)$. Here $G^1_d(\tilde C)$ denotes the Hilbert scheme
of lines in $\tilde C^{(d)}$, that is of curves $P$ such that $P \cdot (x+C^{(d-1)}) =
1$ and $p_a(P)=0$ (cfr. [ACGH] IV.3). \par \noindent Each $P \in \Cal G$ is a pencil
of divisors of degree $2g-2$. If $P$ is general
$\pi^{(2g-2)}$ embeds $P$ in $C^{(2g-2)}$ as a conic $B$. Hence $\pi^{(2g-2)}$ induces
a rational map
$$
h: \Cal G \to \Cal B,
$$
sending $P$ to $B$. It is easy to see that $h(P) = h(i^*P)$ and that $\pi^{(2g-2)}/P$
is an embedding iff $P' \neq i^*P'$, $P'$ being the moving part of $P$. Let
$a: \Cal G \to \tilde J$, $b: \Cal B \to J$ be the natural Abel maps, then
$$
\CD
{\Cal G} @>>h> {\Cal B} \\
@VVaV @VVbV \\
{\tilde J} @>>{Nm}> J \\
\endCD
\tag 2.13
$$
is a commutative diagram. We are interested in the scheme
$$
X = a^*Nm^*(o) = h^*b^*(o), \tag 2.14
$$
where $o$ is $\omega_C$ considered as a point of $J$. In other words $X$
is the is the fibre at $\omega_C$ of the natural norm map $G^1_{2g-2}(\tilde C) \to
Pic^{2g-2}(C)$. Since $Nm^*(o) = P \cup P^-$, $X$ is not connected. The underlying set
of $X$ is naturally identified with the family of  pairs
$$
(\tilde L, \Gamma), 
$$
where $\Gamma$ is a $2$-dimensional subspace of $H^0(\tilde L)$ and $Nm \tilde L =
\omega_C$. Note that the image of $a$ is $Sing \tilde \Theta$. Let $\Xi^- = P^-\cdot
Sing \tilde \Theta$, one has 
$$
X = X^+ \cup X^-,
$$
where
$$ 
X^+ = a^*\Xi \quad \text {,} \quad X^- = a^*\Xi^-. 
$$
\proclaim {2.15 THEOREM} Every component of $X$ is at least $g-2$ dimensional.
Moreover $X^+$ and $X^-$ are integral of dimension $g-2$ if $\pi$ is general 
and $g \geq 4$.
\endproclaim
\demo {Proof} The second statement follows from the first one and the theorems of
Bertram and Welters about reducedness, irreducibility and codimension of
Brill-Noether  loci in a general Prym variety, ([B], [W]). The first statement too is
a consequence of the same results.
\enddemo \noindent
We want to point out that the Hilbert scheme of conics of $\mid \omega_C \mid$ is the
fibre
$$
\Cal B_{\omega_C} = b^{-1}(o)
$$ 
of the map $b$. This scheme is related to  $\bold Q^3$ by the birational map
$$
\delta: \Cal B_{\omega_C} \to \bold Q^3, \tag 2.16
$$ 
which associates to a smooth conic $B$ its dual hypersurface $B^* \subset \bold P^+$.
This defines a map $\alpha = \delta \cdot h/X$ which is fundamental in this paper.
In the beginning we want to define $\alpha$ in a slightly different way.  To do this
we consider on $X$ the universal line
$$
p: \Cal P \to X
$$
induced by the Hilbert scheme $\Cal G$. $\Cal P$ is the projectivization of a rank two
vector bundle $\Cal U$. Moreover the fibre product $\Cal P^* \times_p \Cal P^*$
is embedded in  $\bold P(\Cal U \otimes \Cal U)^*$  as the locus of indecomposable
vectors. Over the point $P$ this embedding is precisely the Segre embedding of $\bold
P(\Cal U_P^*) \times \bold P(\Cal U_P^*)$ as a quadric surface. This defines a section
$$
\sigma: X \to \bold PSym^2 (\Cal U \otimes \Cal U), \tag 2.17
$$
which associates to $P$ such a quadric. Then we consider the multiplication map
$$
\mu: \Cal U \otimes \Cal U \to  H^0(\Cal O_{\bold P}(1))\otimes \Cal O_{\tilde C},
$$
sending $s \otimes t$ to $si^*t$. The $2$-symmetric product of $\mu$ induces
a map
$$
m: \bold PSym^2(\Cal U \otimes \Cal U) \to \bold PH^0(\Cal O_{\bold P}(2)).
\tag 2.18
$$
Taking the restriction of $m(\sigma(P))$ to $\bold P^+$ and $\bold P^-$ one defines
two maps:
$$ 
m^+: \bold PSym^2(\Cal U \otimes \Cal U) \to \bold PH^0(\Cal O_{\bold P^+}(2))
$$
and
$$
m^-: \bold PSym^2(\Cal U \otimes \Cal U) \to \bold PH^0(\Cal O_{\bold P^-}(2)).
$$
Let $P$ be the point $(\tilde L, \Gamma)$ of $X$ and let
$$
\tilde Q = m(\sigma(P)). 
$$
We can describe the quadric $\tilde Q$ in a very concrete way. At first we have $\Cal
U_P = \Gamma$. Fix a basis $\lbrace s_1,s_2 \rbrace$ of $\Gamma$ so that $z_{uv} = s_u
\otimes s_v$ is a basis of $\Gamma \otimes
\Gamma$. Then consider on $\bold P$ the linear forms $x_{uv} = m(z_{uv}) = s_ui^*s_v$
and the matrix
$$
A = (x_{uv}). \tag 2.19
$$
Since the image of $\mid \Gamma \mid^* \times \mid \Gamma \mid^*$ by Segre embedding is
$\lbrace det (z_{uv})=0 \rbrace$, it follows 
$$\tilde Q = \lbrace det A = 0
\rbrace.$$
\bf2.20 REMARK \rm The construction of $\tilde Q$ is classical and well known: see
e.g. [AM], [ACGH], [K]. If $\Gamma = H^0(\tilde L)$ then, by Riemann-Kempf theorem, 
$\tilde L$ is a double point of $\tilde \Theta$ and $\tilde Q$ its projectivized
tangent cone. $\tilde Q$ clearly contains $\tilde C$.  
\proclaim {2.21 LEMMA} Let $M$ and $ \tilde c$ respectively be the moving and
the fixed part of $\mid \Gamma \mid$. Assume $\Lambda \subset \tilde Q$ is a 
maximal linear subspace. Then $\Lambda \cdot \tilde C$ $=$ $\tilde c+i^* \tilde c+f$,
where $f \in M \cup i^*M$. In particular $M$ is cut on $\tilde C$ by a ruling of
$\tilde Q$.
\par
\noindent
\endproclaim
\demo {Proof} We have $s_u = \sigma s'_u$ ($u=1,2$), where
$s'_1$, $s'_2$ generate $M$ and $div (\sigma) = \tilde c$. Moreover a system of
equations for $\Lambda$ can be written as $(a,b)A=(0,0)$ or as $A^t(a,b)=(0,0)$,
for some constants $a,b$. Since $\lbrace x_{uv}=0\rbrace \cdot \tilde C$
$=$ $div (s_ui^*s_v)$, it follows $\Lambda \cdot \tilde C$ $=$ $\tilde c +
i^* \tilde c+div(as'_1+bs'_2)$ or $\Lambda \cdot \tilde C$ $=$ $\tilde c + i^* \tilde
c+div i^*(as'_1+bs'_2)$. 
\enddemo \noindent 
Now we want to describe $m^+(\sigma(P))$ and $m^-(\sigma(P))$, that is
the intersections 
$$
q = \tilde Q \cdot \bold P^+ \quad  \text {and} \quad t = \bold P^- \cdot \tilde Q.
\tag 2.22
$$
At first we recall how $i^*$ acts on a linear form $x$: $i^*x=x$ $\Leftrightarrow$ $x$
is zero on $\bold P^-$, $i^*x=-x$ $\Leftrightarrow$ $x$ is zero on $\bold P^+$.
Secondly we observe that the terms in the matrices
$$
(A+^tA) \quad \text {and} \quad (A-^tA)
$$
are respectively invariant and antiinvariant by $i^*$. This follows from
$i^*x_{uv}=x_{vu}$. Then one easily deduces that the restrictions of $\tilde Q$ 
to $\bold P^+$ and $\bold P^-$ are respectively
$$
q = \lbrace det (A+^tA)=0 \rbrace \quad \text {and} \quad t = \lbrace det
(A-^tA)=0\rbrace.
\tag 2.23
$$
Here we have identified, with some abuse, a linear form on $\bold P$ and
its restriction to $\bold P^+$ or to $\bold P^-$. We will do this when no
confusion is possible. Note that
$$
det (A+^tA) =  4x_{11}x_{22}-(x_{12}+x_{21})^2 = 0 \quad \text {and}
\quad det (A-^tA) = (x_{12}-x_{21})^2,
$$
so $q$ is a quadric of rank $\leq 3$ and $t$ is either $\bold P^-$ or a double
hyperplane in it. $m^-\cdot \sigma$ is not defined at $P$ iff $t = \bold P^-$. 
Since $x_{11}$, $x_{22}$ are linearly independent on $\bold P^+$, $q$ has rank $\geq
2$. Hence $m^+ \cdot \sigma$ is a morphism. 
\proclaim {2.24 PROPOSITION} The following conditions are equivalent: \par \noindent
(i)  $t = \bold P^-$, \par \noindent
(ii) the moving part of $\mid \Gamma \mid$ is the pull-back by $\pi$ of
a pencil on $C$. \par \noindent (iii) $q$ contains $C$.
\endproclaim
\demo {Proof} Let $s_u=\sigma s'_u$ ($u=1,2$) as in the proof of (2.21). Then: (i)
$\Leftrightarrow$ $x_{12}-x_{21}=0$ $\Leftrightarrow$ $s_1i^*s_2 = s_2i^*s_1$ 
$\Leftrightarrow$ $s'_u = i^*s'_u$ ($u=1,2$) $\Leftrightarrow$ (ii). 
Moreover, it is easily seen that (i) $\Leftrightarrow$ $\tilde Q \supset \tilde C \cup
C^-$ $\Leftrightarrow$ $\tilde Q \supset \tilde C \cup C$ $\Leftrightarrow$ (iii).
\enddemo \noindent
\proclaim {2.25 LEMMA} $t \cdot C^- = q \cdot C$. In particular
$q$ belongs to the variety $Y$.
 \endproclaim
\demo {Proof} One immediately checks that $4x_{11}x_{22}-(x_{12}+x_{21})^2$ and
$(x_{12}-x_{21})^2$ restrict to the same bicanonical section of $\tilde C$. Then
the result follows from the equalities
$$
\pi^*(t\cdot C^-) = \tilde C \cdot \lbrace (x_{12}-x_{21})^2 = 0 \rbrace
= \tilde C \cdot \lbrace
4x_{11}x_{22}-(x_{12}+x_{21})^2=0\rbrace = \pi^*(q\cdot C).
$$
\enddemo
\noindent 
By the previous lemma the image of $m^+ \cdot \sigma$ is in $Y$. On the other hand the
image of $m^- \cdot \sigma$ is in $\bold P^{-*}$, where $\bold P^{-*}$ is embedded in
the space of quadrics of $\bold P^-$ by the squaring map. From now on we will adopt
the following
\proclaim {2.26 DEFINITION} \par \noindent
- $\alpha: X \to Y$ is the morphism $m^+ \cdot \sigma$. \par \noindent
- $ \beta: X \to \bold P^{-*}$ is the rational map $m^- \cdot \sigma$.
\endproclaim 
\proclaim {2.27 THEOREM} Let $v: \bold P^{-*} \to V$ be the restriction isomorphism
sending
$t$ to $t\cdot C^-$ and let $\lambda$ be the restriction map considered in 2.6.
Then the diagram
$$
\CD X @>{\alpha}>> Y \\ @VV{\beta}V @VV{\lambda/Y}V \\ {\bold P^{-*}} @>>v> V\\
\endCD 
$$
is commutative.
\endproclaim 
\demo {Proof} It follows immediately from lemma 2.25 and the definition of
$m^+$, $m^-$.\enddemo
\noindent Now we can relate the previous constructions to the Gauss map. First we
recall that  the target space of  this map is $\bold P^{-*}$, because the tangent
sheaf of $P$ is
$\Cal O_P \otimes H^-$. Therefore, the Gauss map
$$
\gamma: \Xi \to \bold P^{-*} \tag 2.28 
$$ 
associates to a point $\tilde L \in \Xi - Sing \Xi$ a hyperplane $\gamma(\tilde L)$
in $\bold P^-$. Since $\tilde L$ is a smooth point, the space $\Gamma =
H^0(\tilde L)$ is two-dimensional, (cfr. [ACGH] p.298). Hence $P = (\tilde L, \Gamma)$
is a point of
$X$. In particular $\beta(P)$ is another hyperplane in $\bold P^-$.
The main fact is now Tjurin's geometric description of $\gamma$ which says (cfr.
[T],1.4):
$$
\gamma(\tilde L) = \beta(P). \tag 2.29
$$
\proclaim {2.30 THEOREM} If $X^+$ is irreducible, then it holds $deg \beta/X^+ = deg
\gamma$.
\endproclaim
\demo {Proof} The Abel map $a/X^+: X^+ \to \Xi$ sends $P=(\tilde L, \Gamma)$
to $\tilde L$. Moreover we have $\beta/X^+ = \gamma \cdot a/X^+$ because
$\gamma(\tilde L) = \beta(P)$. $a/X^+$ is clearly invertible on $\Xi-Sing \Xi$.
If $X^+$ is irreducible, then $a/X^+$ is birational and the result follows.
\enddemo \noindent
The statement is no longer true if $X^+$ is reducible. In this case the irreducible
components of $X^+$ which are contracted by the Abel map $a$ contribute to  the
degree of $\beta/X^+$. An example of this situation will appear in section $4$, where
$C$ is trigonal. From now on we fix the following notations:
\proclaim {2.31 DEFINITION} \par \noindent
$b_+(\pi)$ $=$ degree of $\beta$ on $X^+$, \par \noindent
$b_-(\pi)$ $=$ degree of $\beta$ on $X^-$.
\endproclaim \noindent
It is clear that $b_+(\pi)$ is the degree of the Gauss map if $\pi$ is general.
The remaining part of this paper will be devoted to the solution of the following
problem:
$$
\text {\it Compute $b_+(\pi)+b(_-(\pi)$ and $b_+(\pi)-b_-(\pi)$ if
$\pi$ is general.}
$$
\bigskip \noindent \noindent
\bf 3. \bf The maps b$_{+}$ - b$_{-}$ and b$_{+}$ + b$_{-}$ \rm \par \noindent
Let $\Cal R_g$ be the moduli space of $\pi: \tilde C \to C$. With some abuse, we will
denote in the same way $\pi$ and its moduli point in $\Cal R_g$. We have maps
$$ b_+: \Cal R_g \to \bold Z \quad \text {and} \quad b_-: \Cal R_g \to \bold Z,
$$ sending $\pi$ respectively to $b_+(\pi)$ and to $b_-(\pi)$. It is clear that
the values of $b_++b_-$ and of $b_+-b_-$ are constant on some open dense subset.
In this section we compute this constant value for $b_++b_-$. Then we introduce a
method for computing the constant value of $b_+-b_-$. The method will be applied in
the next section to the case where $C$ is trigonal, this will make the computation
effective .
\proclaim {3.1 LEMMA} The diagram
$$
\CD X @>>h> \Cal B_{\omega_C} \\ @VV{\sigma}V @VV{\delta}V\\ {\bold PSym^2(\Cal U
\otimes \Cal U)} @>>{m^+}> {\bold Q^3}\\
\endCD
$$ is commutative. In particular $\alpha = \delta \cdot h$.
\endproclaim
\demo {Proof} Let $P=(\tilde L, \Gamma)$. Assume $q = (m^+ \cdot\sigma)(P)$ and
$q'=(\delta \cdot h)(P)$. By definition of $\delta \cdot h$, the family $\lbrace <
\pi_* \tilde f >,\quad \tilde  f\in \mid \Gamma \mid \rbrace$ is the family of
tangent hyperplanes to $q'$.  On the other hand we have $< \pi_* \tilde  f > = \bold
P^+ \cdot <\tilde f+i^*\tilde f>$ for each $\tilde f$. Keeping our notations as in
2.19, we can assume $\tilde f = div s_1$.  Then $<\tilde f+i^*\tilde f>=\lbrace
x_{11}=0 \rbrace$, which is tangent to $q$. Hence $<\pi_* \tilde f>$ is tangent to
$q$ so that $q = q'$.\enddemo
\noindent Now we want to describe with some detail the fibres of $\alpha$. Let 
$$
q \in Y,
$$
to avoid unnecessary complications \it we will assume that $q$ has rank three. \rm Let
$$
c \in Div(C) 
$$
be the divisor associated to $C \cdot Sing q$, we consider the family  
$$
B_q = \lbrace c+f \in \mid \omega_C \mid / <c+f> \text{is a hyperplane tangent
to $q$} \rbrace. \tag 3.2
$$
We want to point out that $B_q$ is a smooth conic in $\mid \omega_C \mid$ and that
$\delta(B_q) = q$. Moreover, if $q$ contains $C$, it holds $B_q = c + \pi_*M$, where
$M = \pi^*N$ and $N$ is the base-point-free pencil which is induced by the ruling of
$q$.  
\proclaim {3.3 LEMMA} On $\tilde C$ there exists a base-point-free pencil $M$ such
that 
$$
B = c + \pi_*M = c + \pi_*i^*M.
$$
Moreover $M$ is unique modulo the action of $i^*$.
\endproclaim
\demo {Proof} Assume $M$ exists and is generated by $s_1,s_2$. Let $\tilde Q$ $=$
$\lbrace det (x_{uv}) = 0 \rbrace$, ($u,v=1,2$), where $x_{uv} = \sigma s_ui^*s_v$
and $div \sigma = \pi^*c$. $\tilde Q$ is a quadric in
$\bold P$ such that: (i) $\tilde Q \cdot \bold P^+ = q$, (ii) $\tilde Q \cdot \bold
P^-$ has rank $\leq 1$, (iii) the equation of $\tilde Q$ is $i^*$-invariant, (iv)
$\tilde Q \supset \tilde C$. (ii),(iii),(iv) can be proved for $\tilde Q$ exactly as
for the quadric considered in 2.19. In the same way it follows that $ \lbrace
4x_{11}x_{22}-(x_{12}+x_{21})^2=0 \rbrace$ is a cone over $q'$, where $q' =
\tilde Q \cdot \bold P^+$. One can check that a hyperplane $h$ is tangent to this
cone iff $h = <\pi^*c+\tilde m+i^*\tilde m>$, for some $\tilde m \in M$. On the other
hand, by the proof of 3.1, $<\pi^*c+\tilde m+i^*\tilde m>$ is tangent to $q$.
Then $q' = \delta(B_q) = q$ and (i) holds.  A quadric like $\tilde Q$ satisfies the
following property too: the pencil $\pi^*c+M$ is cut on $\tilde C$ by a ruling of
maximal linear subspaces of $\tilde Q$, the other ruling cuts $\pi^*c+i^*M$.
Therefore $M$ is uniquely defined by $\tilde Q$ up to the action of $i^*$. The
statement then follows if there exists a unique quadric $\tilde Q$
satisfying (i),(ii),(iii),(iv). \par \noindent  Uniqueness: by (iii) $\tilde Q$ has
an equation of type $A+B$, with $A$ $\in$ $H^0(\Cal O_{\bold P^+}(2))$ and $B \in
H^0(\Cal O_{\bold P^-}(2))$. Let $A_1+B_1$ and $A_2+B_2$ be the equations of two
quadrics satisfying (i),(ii),(iii),(iv). By (i) we can assume $A_1=A_2$. By (ii) and
(iv) $B_1-B_2$ has rank $\leq 2$ and is zero on $\tilde C$. This implies $B_1=B_2$. 
\par \noindent
Existence: $\tilde Q$ is the cone of vertex $\bold P^-$ over
$q$ if $q \supset C$. If not $q \cdot C = 2d$ and $d$ is cut on $C^-$ by a hyperplane.
Let $B$ be a squared equation of it, $A$ an equation of $q$. Regarding $A$, $B$ as in
$H^0(\Cal O_{\bold P}(2))$, we consider the pencil  $S = \lbrace z_0A+z_1B = 0
\rbrace$. Since
$\lbrace A=0 \rbrace \cdot \tilde C = \pi^*2d = \lbrace B=0 \rbrace
\cdot \tilde C$, a quadric of $S$ contains $\tilde C$. This is $\tilde Q$.
\enddemo 
\proclaim {3.4 COROLLARY} A point $P=(\tilde L, \Gamma)$ belongs to $\alpha^{-1}(q)$
if and only if $\mid \Gamma \mid$ is $\tilde c + M$ or $\tilde c +i^*M$, for some
divisor $\tilde c$ on $\tilde C$ satisfying $\pi_*\tilde c = c$.
\endproclaim
\demo {Proof} Immediate from the lemma. \enddemo
\proclaim {3.5 PROPOSITION} Let $q \in Y$ and let $n = deg c$. Then: \par \noindent
(1) $\alpha^*(q)$ is a zero-dimensional scheme of length $2^{n+1}$.
\par \noindent
(2)  Assume $c>0$. Then $\alpha^*(q)\cdot X^+$ and $\alpha^*(q)\cdot X^-$ have the
same length.
\endproclaim
\demo {Proof} Let $B = \delta^*(q)$. $B$ is a point of the Hilbert scheme $\Cal
B_{\omega_C}$, that is a smooth conic in $\mid \omega_C \mid$. We consider $Z =
\pi^{(2g-2)*}(o_B)$, where $o_B$ is the generic point of $B$ and $\pi^{(2g-2)}: \tilde
C^{(2g-2)} \to C^{(2g-2)}$ is the push-down map. Since $h: \Cal G \to \Cal B$ is the
map of Hilbert schemes induced by $\pi^{(2g-2)}$, it follows that $(\delta
\cdot h)^*(q)$ is a subscheme of $Z$. The support of $(\delta \cdot
h)^*(q)$ is the set of generic points $o_P$ of lines $P$ such that 
$\pi^{(2g-2)}(P) = B$. The multiplicity of $o_P$ in $(\delta \cdot h)^*(q)$ is 
the ramification index $\nu_P$ of $\pi^{(2g-2)}$ at $o_P$. Therefore $\nu_P$
is the index of $\pi^{(2g-2)}$ at a general point of $P$ times the degree of
$\pi^{(2g-2)}/P$. To compute $\nu_P$ we consider the base extension
$$
\CD
{\tilde C^{(n)} \times \tilde C^{(m)}} @>>> {\tilde C^{(2g-2)}} \\ 
@V{\pi^{(n)}\times \pi^{(m)}}VV @V{\pi^{(2g-2)}}VV \\
{C^{(n)} \times C^{(m)}} @>>> {C^{(2g-2)}} \\
\endCD
$$
where the horizontal arrows are the sum maps and $m+n=2g-2$. By lemma 3.3 $o_P$
is in $(\delta \cdot h)^*(q)$ iff $P = \tilde c + M$ or $P = \tilde c + i^*M$,
with $M$ base-point-free and $\pi_* \tilde c = c$. From the
equality $P = \tilde c + M$ and the previous diagram it follows that $\nu_P$ is  the
ramification index of $\pi^{(n)}\times \pi^{(m)}$ at the generic point of $\lbrace
\tilde c \rbrace \times M$. Then $\nu_P = \nu_{\tilde c} \nu_M$, where $\nu_{\tilde
c}$ is the index of $\pi^{(n)}$ at $\tilde c$ and  $\nu_M$ is the index of $\pi^{(m)}$
at the generic point of $M$. Since $M$ is base-point-free, $\pi^{(m)}$ is not
ramified along $M$. Hence $\nu_M$ coincides with the degree of $\pi^{(m)}/M$ onto  its
image. This is one if $M \neq i^*M$, two if $M=i^*M$.  We can conclude that:
length$(\delta \cdot h)^*(q)$ $=$ $\Sigma \nu_P$ $=$
$\Sigma \nu_{\tilde c} \nu_M$ $=$ $2\Sigma\nu_{\tilde c}$ $=$ $2 deg \pi^{(n)}$ $=$
$2^{n+1}$. This shows (1), since $\alpha = \delta \cdot h$. \par \noindent 
(2) $\pi^{(n)*}(c)$ embeds in $\alpha^*(q)$ via the map sending $\tilde c$ to $\tilde
c+M$. Let $U$ be the image of such an embedding. If $M \neq i^*M$,
$\alpha^*(q)$ is the disjoint union of $U$ and $i^*U$. If $M=i^*M$,
$\alpha^*(q)$ is $U$ with a double structure, because $\nu_M=2$. On the other hand we
have $\pi^{(n)*}(c) \subset S$, where $S = \pi^{(n)*}\mid c \mid.$ $S$ is the disjoint
union of $S^+$ and $S^-$, where $S^+ = \lbrace \tilde c \in S/ \tilde c + M \in X^+
\rbrace $ and $S^-=S-S^+$. If $c$ is the zero divisor $S$ is trivially one point and
$ \lbrace S^+,S^- \rbrace = \lbrace S,\emptyset \rbrace$. If $c > 0$, $S^+$ and $S^-$
are well known special varieties. They have been studied by Beauville in [B2]. In
particular they are numerically equivalent in $\tilde C^{(n)}$, ([B2] section 1).
Hence $\pi^{(n)*}(c) \cdot S^+$ and $\pi^{(n)*}(c) \cdot S^-$ have the same
length. This implies (2).
\enddemo \noindent
\proclaim {3.6 COROLLARY} $\alpha$ is a finite morphism onto $Y$. \endproclaim
\demo {Proof} It follows from statement 3.5(1). \enddemo
\proclaim {3.7 PROPOSITION} Assume $\pi$ is standard and Petri general. Then the
degree of $\alpha$ is two.
\endproclaim 
\demo {Proof} By assumption no element of $Y$ contains $C$, moreover $X$ and
$Y$ are reduced of dimension $g-2$. Let $S = \lambda^*(d)$ be a general
fibre of $\lambda/Y$, then $S$ is smooth and $dim S = 0$.  $S$ is the subset of rank
three quadrics in the linear system $\mid I_{2d}(2) \mid$ of all quadrics cutting $2d$
on $C$.  By proposition 3.5(1) it suffices to show that 
$Sing q \cap  C = \emptyset$, for some $q \in S$. Let $T$ be the embedded 
tangent space to $S$ at $q$ and let $\Cal I_{Sing q}$ be the Ideal of
$Sing q$, then $T = \mid I_{2d}(2) \mid \cap \mid \Cal I_{Sing q}(2) \mid$,
(cfr.[ACGH] A-5 p.101). In other words $T = \bold PKer \rho$, where $\rho:H^0(\Cal
I_{2d}(2)) \to H^0(\Cal O_{Sing q}(2))$ is the restriction map. Note that
$h^0(\Cal I_{2d}(2))-h^0(\Cal O_{Sing q}(2)) = 1$. Assume  $x \in Sing q \cap C$,
then $\rho$ is not surjective and hence $dim T \geq 1$: against the smoothness of $S$.
\enddemo \noindent
\proclaim {3.8 THEOREM } $b_+(\pi)+b_-(\pi) = 2D(g)$ if $\pi$ is standard and Petri
general.
\endproclaim
\demo {Proof} By assumption $X$, $Y$ are reduced. Moreover $deg \alpha = 2$ and $deg
\lambda/Y = D(g)$ so that $deg (\lambda \cdot\alpha) = 2D(g)$. From $\beta = \lambda
\cdot \alpha$ one has $b_+(\pi)+b_-(\pi) = deg \beta = 2D(g)$.
\enddemo \noindent
Now we study the difference map $b_+-b_-$. We will work under the following
\proclaim {3.9 ASSUMPTION} \par \noindent
(1) Let $C \subset q$, where $q$ is a rank three quadric. Then $Sing q \cap C \neq
\emptyset$. 
\par \noindent
(2) No irreducible component of $Y$ is in the center of $\lambda$.
\par \noindent
(3) The general fibre of $\lambda/Y: Y \to V$ is finite.
\endproclaim \noindent
Condition (1) guarantees that $\alpha^*(q) \cdot X^+$ and $\alpha^*(q)\cdot X^-$ have
the same length as soon as $q$ contains $C$. (1) is not satisfied iff there exists a
globally generated theta characteristic $\theta$ on $C$. (2) and (3) seem to be
always true.  \par \noindent
\proclaim {3.10 DEFINITION} $R \subset \Cal R_g$ is the subset where the
assumption holds.
\endproclaim \noindent
Note that (1),(2),(3) are satisfied if no rank three quadric contains $C$. Hence $R$
contains $R_{Petri}$. In particular \it $R$ contains a dense open subset. \rm
A good reason for considering $R$ instead of $R_{Petri}$ is that $R$ has non
empty intersection with the locus of points $\pi$ such that $C$ is trigonal:
see section 4. Now we fix a smooth family
$$
\pi: \tilde \Cal C \to \Cal C \tag 3.11
$$
of \'etale double coverings of genus $g$ curves. This means that
$\pi$ is a connected \'etale double covering of a smooth family 
$$
p: \Cal C \to T \tag 3.12
$$
of irreducible genus $g$ curves, we will assume that each fibre of $p$ is not
hyperelliptic. $\pi$ is defined by $\eta \in Pic$ $\Cal C$, $T$ is a smooth integral
affine curve. For each $t \in T$, $\pi$ induces an \'etale double covering
$$
\pi_t: \tilde \Cal C_t \to \Cal C_t. 
$$ 
Using the previous data we construct the families of symmetric products
$$
\tilde p^{(n)}:\tilde \Cal S^{(n)} \to T \quad \text {and} \quad p^{(n)}:  \Cal S^{(n)}
\to T,
$$
with fibres  $\tilde \Cal C_t^{(n)}$ and $\Cal C_t^{(n)}$. Then we consider the
Hilbert scheme $\Cal B$ of conics in $\Cal S^{(2g-2)}$ and the Hilbert scheme  
$\Cal B_{\Cal H}$ of conics in $\Cal H$, where $\Cal H = \bold Pp_*\omega_{\Cal
C/T}$. Due to the equality  $\Cal H_t = \mid \omega_{\Cal C_t} \mid$ there is a
natural inclusion $\Cal H \subset \Cal S^{(2g-2)}$, so that $\Cal B_{\Cal H} \subset
\Cal B$. Finally we also consider the Hilbert scheme $\Cal G$ of lines in $\tilde \Cal
S^{(2g-2)}$. All these are schemes over $T$, they admit standard definitions we
have omitted for brevity. Let 
$$
\phi: \tilde \Cal S^{(2g-2)} \to \Cal S^{(2g-2)}
$$ 
be the morphism sending $d \in \tilde \Cal C_t^{(n)}$ to $\pi_{t*}d \in \Cal
C^{(n)}_t$. Then $\phi_*$ induces a map
$$
h: \Cal G \to \Cal B
$$
sending a general line $P \in \Cal G$ to the conic $h(P) = \phi_*P$. By
definition we put  
$$
\Cal X = h^* \Cal B_{\Cal H}. \tag 3.13
$$
For each $t$ $\Cal X_t$ is our usual variety $X$, defined from $\pi_t$ as in 2.14. We
have 
$$
\Cal X = \Cal X^+ \cup \Cal X^-, \tag 3.14
$$ 
$\Cal X^+$ and $\Cal X^-$ being disjoint components. Each component of $\Cal X$ is at
least $g-1$ dimensional: the same proof used in 2.15 works. By definition $P$ is a
point of $\Cal X^+$  iff $P$ is a line in a complete linear system of odd dimension. 
Let 
$$
\Cal Q = \bold PSym^2 p_*\omega_{\Cal C/T}
$$
and let $\Cal Q^3 \subset \Cal Q$ be the locus of quadrics of rank $\leq 3$. As in
2.16 we define a map
$$
\delta: \Cal B_{\Cal H} \to \Cal Q,
$$
sending a smooth conic to its dual hypersurface. Then we define
$$
\alpha = \delta \cdot h/ \Cal X : \Cal X \to \Cal Q. \tag 3.15
$$
$\Cal Q^3$ is the birational image of $\delta$. Let us fix a
general very ample section
$$ 
s: T \to \bold Pp_*(\omega_{\Cal C/T} \otimes \eta).
$$
Squaring $s$ we obtain the family $\lbrace 2d_t \in \mid \omega_{\Cal C_t}^{\otimes
2} \mid, d_t = s(t) \rbrace.$ Then we construct, in the usual standard way, a vector
bundle
$$\Cal F \subset Sym^2 p_* \omega_{\Cal C/T}\tag 3.16 $$  such that $\Cal F_t \subset
Sym^2H^0(\omega_{\Cal C_t})$ is the space of quadratic forms vanishing on $2d_t$.
Here, and in what follows, each $\Cal C_t$ is canonically embedded. Let
$$
F = \bold P\Cal F,
$$
we remark that each $F_t$ contains a hyperplane $I_t$ which is the linear system of
quadrics through $\Cal C_t$. $I_t$ is the fibre at $t$ of a projective bundle $I
\subset \Cal Q$ and $I$ is the indeterminacy locus of the natural restriction map
$$
\lambda: \Cal Q \to \bold Pp_*\omega_{\Cal C/T}^{\otimes
2}. \tag 3.17
$$
The resolution of the indeterminacy is $\sigma \cdot \lambda$, $\sigma: \Cal Q' \to
\Cal Q$ being the blowing up of $I$. Let $\Cal Y'$ be the strict transform
of $\Cal Y = \alpha( \Cal X)$, we consider the commutative diagram
$$
\CD
{\Cal X^{'+} \cup \Cal X^{'-}} @>{\alpha'}>> {\Cal Y'} @>{\lambda'/\Cal Y'}>>
{\bold Pp_*\omega_{\Cal C/T}^{\otimes 2}}\\
@VVV @V{\sigma/\Cal Y'}VV @V{id}VV\\
{\Cal X^+\cup \Cal X^-} @>>{\alpha}> {\Cal Y} @>>{\lambda/\Cal Y}> {\bold
Pp_*\omega_{\Cal C/T}^{\otimes 2}}\\
\endCD
$$
Here $\lambda': \Cal Q' \to \bold Pp_*\omega_{\Cal C/T}^{\otimes 2}$ is the 
projective bundle structure induced by $\lambda$ and the left square is a
base extension. Recall that each component of $\Cal X$ has  dimension $\geq g-1$,
the same  property holds for $\Cal X^{'+} \cup \Cal
X^{'-}$ because $\alpha$ is finite and $\sigma$ is birational. Let 
$$
F' \subset \Cal Q'
$$
be the strict transform of $F$ by $\sigma$. Since $I$ is a divisor in $F$, 
$\sigma/F': F' \to F$ is a  biregular map. In particular $F'$ is locally a 
complete intersection in $\Cal Q'$ and its codimension is $g-2$. This implies 
that each irreducible component of
$$
S = {\alpha'}^*F'.
$$ 
has dimension $\geq 1$. Finally we fix a point $o \in T$ and make the
following \it $\underline{assumption}$: \rm
\bigskip \noindent 
(i) $\pi_t$ is Petri general if $t \neq 0$, so that no rank $3$ quadric contains $\Cal
C_t$,  \par \noindent   
(ii) $d_t$ is general in its Prym-canonical linear system, in
particular $\lambda^*(d_o)$ is finite, \par \noindent   
(iii) $\pi_o$ is in the subset $R \subset \Cal R_g$, defined in 3.10.
\proclaim {3.18 LEMMA} Each irreducible component of $S$ is a curve which maps
onto $T$ via the rational map $\lambda' \cdot \alpha'$. \endproclaim
\demo {Proof} It suffices to show that $S_t = {\alpha'}^*F'_t$ is finite.
Let $t\neq o$, one can check in the previous diagram that $\alpha^*F_t$ is
biregular to $S_t$. By (i) $I_t \cap \Cal Q^3 = \emptyset$, hence  $F_t \cdot \Cal
Q^3$ is finite as well as $\alpha^*F_t$. Let $t=o$, consider the strict transform
$\Cal Y"$ of $\Cal Y_o$ and the Zariski closure  $E$ of $\Cal Y'_o-\Cal Y"$. Then
it holds $S_o = {\alpha'}^*(F'_o \cdot \Cal Y") \cup {\alpha'}^*(F'_o \cdot
E)$. The scheme $F'_o \cdot \Cal Y"$ is the fibre over the point $2d_o$ of the map
$\lambda/ \Cal Y_o: \Cal Y_o \to \mid \omega_{\Cal C_o}^{\otimes 2} \mid$. Therefore
it is finite because,  by (ii) and (iii), $d_o$ is general and $\pi_o$ is in $R$. The
same holds for ${\alpha'}^*(F'_o
\cdot \Cal Y'_o)$ because $\alpha'$ is a finite map. Assume $E$ is not empty and
consider the Zariski closure $\Cal Z$ of  $\cup \Cal Y'_t$, $t \neq o$. Then $E$ is 
a component of the fibre $\Cal Z_o$. By (i) $\Cal Z_t$ is split in two distinguished
irreducible components of dimension $g-2$, for each $t \neq o$. This implies $dim E =
g-2$. Now observe that $F'_o$ is a general fibre of a projective bundle over the
$g-2$-dimensional projective space $\mid \omega_{\Cal C_o} \otimes \eta \mid$. Since
$dim E = g-2$ it follows that ${\alpha'}^*(F'_o \cdot E) $ is finite. This completes
the proof.
\enddemo \noindent 
Let 
$$
S^+ = S \cdot\Cal X^{'+} \quad \text {and} \quad S^- = S \cdot\Cal X^{'-},
$$
we consider on $F'$ the family of zero-cycles
$$
\delta_t =  [{\alpha}^{'}_*S^+][F'_t]-[{\alpha}^{'}_*S^-][F'_t], \quad t \in T.
\tag 3.19
$$
The degree of $\delta_t$ is constant, moreover $deg \delta_t = b_+(\pi_t)-b_-(\pi_t)$
if $t \neq o$. On the other hand we can write
$$
\delta_o = e + \delta'_o,
$$
where $e$ is supported on points $q'$ such that $\sigma(q') \in I_o$ and
$Supp (e) \cap Supp (\delta'_o) = \emptyset$.
\proclaim {3.20 LEMMA} $Supp (e)$ is empty, so that $e = 0$ \endproclaim
\demo {Proof} Let $q = \sigma(q')$, where $q \in I_o$. Then ${\alpha}^{'*}(q')$ and
$\alpha^*(q)$ are biregular: this just follows from the definition of base-extension.
Moreover $S^+\cdot \alpha^{'*}(q')$  and $S^- \cdot \alpha^{*'}(q')$ have the same
length. This follows because, by proposition 3.5(1),  $\Cal X^+ \cdot \alpha^*(q)$ and
$\Cal X^- \cdot \alpha^*(q)$ have the same length. Then $q'$ is not in
$Supp (e)$, which must be empty.
\enddemo \noindent 
We have shown that degree of $\delta'_o$ is exactly $b_+(\pi_o)-b_-(\pi_o)$,
therefore
$$
deg \delta_t = b_+(\pi_t)-b_-(\pi_t), \quad \forall t \in T.
$$
In particular $b_+-b_-$ is \it constant \rm on the curve $T' = \lbrace
\pi_t, t\in T \rbrace \subset R$. Let $\rho$ be any point of  $R$, it is clear
that there are families $\pi: \tilde \Cal C \to \Cal C$ satisfying the previous
assumption and moreover such that $\rho = \pi_o$. Then any two points are connected
by a chain of curves as $T'$ and the next theorem follows.
\proclaim {3.21 THEOREM} The map $b^+-b^-$ is constant on $R$.
\endproclaim \noindent
The constant value of $b^+-b^-$ on $R$ is $2^{g-2}$: this is the argument
of the next section.
\bigskip \noindent
\bf 4. Trigonal Curves and Hyperelliptic Thetas.\rm
\par \noindent  In this section we study the case where $\pi: \tilde C
\to C$ is an \'etale double covering of a \it general trigonal \rm curve. 
We will see that in this special case the fibre of the map  $\lambda: Y \to V$ is
naturally related to the set of theta-characteristics on a hyperelliptic curve of
genus $p = g-3$. Using this relation we will be able to compute that
$b^+(\pi)-b^-(\pi) = 2^{g-2}$ and to show that this is the value of $b^+-b^-$ at a
general point of $\Cal R_g$. To simplify the exposition we will assume $g \geq 4$,
leaving the case $g = 3$ as an exercise. Since $C$ is trigonal~one~has
$$
C \subset R. \tag 4.1
$$
$R$ is a \it general rational normal scroll \rm in the canonical space of $C$, that
is  the \it Hirzebruch surface $F_0$ ($g$ even) or $F_1$ ($g$ odd).\rm The Picard group
of
$R$ is $\bold Z [h] \oplus \bold Z [f]$, where $h$ is a hyperplane section and $f$ is
a fibre of $R$. One computes 
$$ C \sim 3h - (g-4)f \quad , \quad K_R \cong -2h +(g-4)f.$$ \bigskip \noindent
Let $I_R$, ($I_C$), be the linear system of quadrics through $R$,
(through $C$), then
$$ I_R = I_C. \tag 4.2$$ 
In particular $I_R$ is the center of the linear projection $\lambda: \bold Q \to
\bold B$ defined in (2.5). As we will see $I_R$ intersects $Y$, so that the degree
of $\lambda/Y$ is different from $D(g)$. 
\bigskip \noindent
\proclaim {4.3 LEMMA} (i) Let $q \in Y \cap I_C$, then $Sing q \cap C$ is not
empty. \par \noindent
(ii) The dimension of $Y \cap I_C$ is $g-5$.
\endproclaim
\demo {Proof} (i) Since each line in $q$ intersects $Sing q$, the same holds for each
fibre $f$ of $R$. This implies that $Sing q \cap R$ contains a curve $b$. On the other
hand we have $C \sim h-K_R$, which is easily seen to be ample on $R$. Hence $C \cap b
\neq \emptyset$. \par \noindent (ii) Restricting to $R$ the hyperplanes through
$Sing q$ we obtain a net of divisors 
$N \subset \mid h \mid$. Let $b$ be the fixed curve of $N$, it is  easy to see that
$bf=1$. Then $h-b \sim mf$ and the moving part of $N$ is a base-point-free net $M$ in
$\mid mf \mid$. Let $\phi_M: R \to \bold P^2$ be the morphism defined by $M$. We
observe that $\phi_M$ is simply the restriction to $R$ of the linear projection of
center $Sing q$. This implies that $\phi_M(R)$ is a smooth conic. But then $\phi_M =
v \cdot \phi_P$, where $v$ is the $2$-Veronese embedding of
$\bold P^1$ and $\phi_P$ is defined by a pencil $P \subset \mid kf \mid$, $m=2k$.
$(P,b)$ is a point of $G_k \times \mid h-2kf\mid$, where $G_k$ is the Grassmannian
of pencils in $\mid kf \mid$. Inverting the construction, one can easily see that any
pair $(P,b) \in \cup_k G_k \times \mid h-2kf \mid$ defines a rank three quadric
containing $R$ and hence $C$. This yelds a surjective 
$$
f: \cup_k (G_k \times \mid h-2kf \mid) \to (I_R \cap Y)
$$
sending $(P,b)$ to $q$. Counting dimensions it follows $dim G_k \times \mid h-2kf
\mid$ $\leq g-5$.
\enddemo
\noindent Now we consider the restriction map  $r:\mid 2h \mid \to \mid
\omega_C^{\otimes 2} \mid$. Since $2h-C \sim h+K_R$, we have $h^0(\Cal O_R(2h-C)) =
h^1(\Cal O_R(2h-C)) = 0$. Then, from the long exact sequence~of
$$
0 \to \Cal O_R(2h-C) \to \Cal O_R(2h) \to \Cal O_C(2h) \to 0,
$$
it follows that $r$ is an isomorphism. We are interested in the family of curves
$$
V_R = r^{-1}(V), \tag 4.4
$$ where $V$ is the image of  $\mid \omega_C \otimes \eta \mid$ under the 'squaring
map'. An element $E \in V_R$ is a hyperelliptic curve of arithmetic genus $p = g-3$.
$E$ satisfies the condition $E \cdot C = 2d$, for a given  $2d \in V$. The
hyperelliptic map of $E$ is induced by the natural projection of $R$ onto
$\bold P^1$. \par \noindent
\proclaim {4.5 LEMMA} The general $E \in V_R$ is a smooth irreducible curve.
\endproclaim
\demo {Proof} It suffices to produce a smooth $E_0 \in \mid 2h \mid$
and a stable canonical curve $C_0 \sim C$ such that: (i) $E_0 \cdot C_0 = 2d_0$,
(ii)
$Supp d_0
\cap Sing C_0 = \emptyset$, (iii) $ dim \mid d_0 \mid = g-2$. Indeed $(C_0,d_0)$
deforms in a family of pairs $\lbrace (C_t,d_t),t \in T\rbrace$, such that $C_t$ is
smooth and $d_t$ is a Prym-canonical divisor. Moreover for each $d_t$ there exists
exactly one $E_t \sim E_0$ such that $E_t\cdot C_t = 2d_t$. The general $E_t$ is
smooth because $E_0$ is smooth: this implies the lemma. For brevity we only sketch
the construction of a pair $(C_0,d_0)$. Fix $E_0 = q \cdot R$, where $q$ is a rank
three quadric transversal to $R$. A general tangent hyperplane to $q$ cuts on $R$ a
smooth curve $A_0$ such that $A_0 \cdot E_0 = 2a_0$. In the very ample linear system
$\mid C-A_0 \mid$ choose a smooth $B_0$ which is transversal to $A_0$
and moreover satisfies $B_0 \cdot E_0 = 2b_0$, $B_0 \cap A_0 \cap E_0
= \emptyset$. The required pair is $(C_0,d_0)$, with $C_0=A_0+B_0$ and $d_0 =
a_0+b_0$. Note that $C_0$ is the union of a rational and an elliptic curve
intersecting transversally in $g$ points.
\enddemo \noindent
Let $I_{2d}$ be the linear system of quadrics through $2d = E \cdot C$ and let 
$I_E$ be the linear system of quadrics through $E$. We want to point out that
$$ 
I_{2d} = I_E. \tag 4.6
$$
This follows from $I_E \subseteq I_{2d}$ and the fact that $dim I_E = dim I_{2d} =
dim I_R+1$, because $I_R=I_C$. We observe that the fibre at a  general $2d$ of 
the map
$\lambda/Y: Y \to V$ is
$$
\hat Q_E = (I_E-I_R) \cdot \bold Q^3.
$$ 
\proclaim {4.7 LEMMA} \it Let $T_q$ be the projective tangent space to $q \in \hat
Q_E$. If $dim T_q > 0$ then $Sing q \cap R \neq \emptyset$.  \endproclaim
\demo {Proof}  The statement is obvious if $q$ has rank $\leq 2$, therefore we assume
$rk(q)=3$. The embedded tangent space to $\hat Q_E$ at $q$ is the linear system
$T_q$ of all quadrics through $E \cup Sing q$. If $dim T_q > 0$ it follows
$T_q \cap I_R \neq \emptyset$. Hence there exists a quadric $q'$ of rank $\leq 6$
which contains both $R$ and $Sing q$. We claim that this implies $Sing q \cap R \neq
\emptyset$. The proof is immediate if $q'$ has rank $\leq 4$. Assume $q'$ has rank
$6$ and consider the projection $\sigma: q' \to \bold P^5$ of center $Sing q'$ and
the smooth quadric $G = \sigma(q')$. Note that $Sing q' \subset Sing q$ and that
$\sigma (Sing q)$ is a plane. Moreover it can be assumed that $\sigma(R)$ is a surface
too: if not $R \cap Sing q$ is already not empty. $\sigma(R)$ is not a plane because
$R$ is not degenerate, hence $\sigma(R)$ has positive intersection index with $\sigma
(Sing q)$ and this implies $Sing q \cap R \neq \emptyset$. If $q'$ has rank $5$ the 
proof is analogous.
\enddemo 
\proclaim {4.8 COROLLARY} $\hat Q_E$ is smooth and finite if $E$ is
smooth. In particular, every irreducible component $Z \subset Y$ 
such that $\lambda(Z) = V$ is reduced.
\endproclaim
\proclaim {4.9 PROPOSITION} Assume $C$ is a general trigonal curve then $\pi$
belongs to the set $R$ defined in 3.10.
\endproclaim 
\demo {Proof} By the previous lemmas 4.8 and 4.5 the general fibre of
$\lambda/Y$ is finite so that condition 3.9 (3) holds for $\pi$.
By lemma 2.11 and 4.4 conditions 3.9 (1) and 3.9(2) are satisfied too. Hence
$\pi$ belongs to $R$. \enddemo \noindent
Now we want to compute the degree of $\lambda$ on the two components $Y^+ =
\alpha(X^+)$ and $Y^- = \alpha(X^-)$ of $Y$. Let $E$ be a smooth element
of $V_R$, $E$ is a hyperelliptic curve. We will show that there is a natural
identification between the set of even (odd) theta's on $E$  and the fibre of
$\lambda/Y^+$ (of $\lambda/Y^-$). We will also see that both $Y^+$ and $Y^-$ are
reducible. Indeed $Y$ contains the components $S_k$, where $S_k$ is defined by the
following condition: the fibre of $\lambda/S_k$ over any smooth $E$ is the subset
of $\theta$'s satisfying $h^0(\theta) = k$. \par \noindent We associate to $E$ the
following sets: 
\par \noindent (4.10) \it the set $T_E$ of theta characteristics on $E$, \rm \par
\noindent   (4.11)  \it the set $Q_E$ of line bundles $\Cal A \in Pic^{p+1}(E)$ such
that
$\Cal A^{\otimes 2} \cong \Cal O_E(h)$, \rm \par \noindent  (4.12) \it the fibre $\hat
Q_E$ of $\lambda/Y$ at $E$. \rm \par \noindent The condition $E \cdot C = 2d$ defines
a bijection
$$
\epsilon_1: Q_E \to T_E, \tag 4.13
$$ 
such that $\epsilon_1(\Cal A) = \Cal A^{\otimes 3}(-d)$. Indeed it holds
$(\epsilon_1(\Cal A))^{\otimes 2} \cong \Cal O_E(3h-C) \cong \omega_E$. 

\proclaim {4.14 LEMMA} Let $C$ be general and let $E$ be general in $V_R$.
Then $E$ satisfies the following property: (*) each $\Cal A \in Q_E$ is globally
generated and has $h^0(\Cal A) = 2$.
\endproclaim 
\demo {Proof} It suffices to prove the statement for a general smooth 
$E \in \mid 2h \mid$, because such an $E$ is in some $V_R$ as above. 
In the Hilbert scheme of $E$ we have the subscheme $\Cal E$ parametrizing non
degenerate curves which are biregular to $E$. One can easily show that $\Cal E$
contains a dense open set $U_1$ of elements $E'$ satisfying (*). On the other hand 
let $\Cal R$ be the Hilbert scheme of $R$. We consider the morphism  $\psi:\Cal E \to
\Cal R$ which is so defined: $\psi(E')=R'$ $=$  $\cup_{x \in E'} <x,j(x)>$, where $j$
is the hyperelliptic involution of $E'$. $R'$  is a Hirzebruch surface $F_n$. We
recall that the locus of points $R' \in \Cal R$ for which $n$ is minimal, (i.e. $n=0$
or $1$), is an open set $U_2$. An element of $U_1 \cap \psi^{-1}(U_2)$ is
projective-isomorphic to some $E' \in \mid 2h \mid$ satisfying (*).
\enddemo \noindent 
In the following \it we will assume that $E$ is smooth and general
as above. \rm Let $q \in \hat Q_E$. Restricting the ruling of $q$ to $E$ we
obtain a base-point-free pencil of divisors of degree $p+1$. This yelds a line bundle
$\Cal A \in Q_E$. By definition
$$
\epsilon_2: \hat Q_E \to Q_E \tag 4.15
$$ 
is the injective map sending $q$ to $\Cal A$. Let $\Cal A \in Q_E$, consider the
rank $3$ quadric
$$ 
q = \cup <a>, \quad a \in \mid \Cal A \mid. 
$$ 
Since $\mid \Cal A \mid$ is base-point-free $Sing q \cap E$ is empty. Hence
$q$ does not contain $R$ and moreover belongs to $\hat Q_E$. Then $\epsilon_2(q) =
\Cal A$ and $\epsilon_2$ is bijective. 
\proclaim  {4.16 PROPOSITION} The degree of $\lambda/Y: Y \to V$ is $2^{2p}$. 
\endproclaim
\demo {Proof}  $\hat Q_E$ is the fibre of $\lambda/Y$. $\epsilon_1 \cdot
\epsilon_2: \hat Q_E \to T_E$ is bijective. $\#T_E = 2^{2p}$. \enddemo
\noindent Now we want to study the theta characteristic
$$
\theta_q = (\epsilon_1 \cdot \epsilon_2)(q) \tag 4.17
$$ for any $q \in \hat Q_E$. Let 
$$
\lbrace (\tilde L_q, \Gamma) , (i^*\tilde L_q, i^*\Gamma) \rbrace = \alpha^*(q),
\tag 4.18 
$$ in particular we want to understand the relation between $\theta_q$ and the line
bundle~$\tilde L_q$. \bigskip \noindent
\proclaim {4.19 THEOREM} For each $q \in \hat Q_E$ one has
$$ h^0(E,\theta_q) = 2 + h^0(\tilde C,\tilde L_q) = 2 + h^0(\tilde C, i^*\tilde L_q).
$$
\endproclaim \noindent Before of giving the proof let us deduce
the main consequences of the theorem. At first it follows that $q$ is in the fibre of
$\lambda/Y^+$ if and only if $\theta_q$ is even. Therefore, recalling the number of
even and odd thetas, we obtain 
\proclaim{4.20 COROLLARY} $deg \lambda/Y^+ = 2^{p-1}(2^p+1)$, $deg \lambda/Y^-
= 2^{p-1}(2^p-1)$. \endproclaim \noindent
Secondly, since $E$ is hyperelliptic, we have a partition of
the set of its thetas according to the dimension of the space of global sections.
This defines a splitting of $Y$ in more than two components. Let
$$ S_k = \overline {\lbrace q \in Y /\text { $q \in \hat Q_E$ for a smooth $E$,
$h^0(\theta_q) = k$} \rbrace}
\tag 4.21
$$ where the overline denotes Zariski closure in $Y$, from the description of theta
characteristics on a hyperelliptic curve it follows
$$
\#(S_k \cap \hat Q_E) = \binom {2p+2}{p+1-2k}, \quad k = 0, \dots [\frac
{p+1}2].
$$ Each $S_k$ is a component of $Y$ and they are reduced by lemma 4.8.
\par \noindent
\bf 4.22 REMARK \rm The degree of $\lambda/S_k$ is $\binom {2p+2}{p+1-2k}$. 
Applying the theorem, we can also remark that
$$
\alpha^*(S_k) = \overline {\lbrace (\tilde L, \Gamma) \in X / h^0(\tilde L) = 2+k
\rbrace}. \tag 4.23
$$
This implies that $\alpha^*(S_0)$ is birational to the Prym-Theta
divisor and that the degree of its Gauss map is $\binom {2p+2}{p+1}$. It is
not surprising that this is also the degree of the Gauss map for the theta
divisor of a $g-1$-dimensional Jacobian. Since $C$ is trigonal, we
know that the Prym of $\pi:\tilde C \to C$ is indeed a Jacobian, ([R]). 
\bigskip \noindent
\demo {PROOF OF THE THEOREM 4.19} Let us fix some preliminary constructions. \par
\noindent (1) \it The Castelnuovo surface of
$E$. \rm
\par \noindent Let 
$$
\rho : \overline R \to R \tag 4.24
$$ be the double covering of $R$ branched on $E$. $\overline R$ is a smooth rational
surface endowed with the pencil of rational curves $\mid \overline f \mid = \mid
\rho^*f \mid$. One can show that 
$$
\mid \overline h \mid = \mid \rho^*h \mid \tag 4.25
$$ 
defines a natural embedding
$$
\phi: \overline R \to \bold P^g, \tag 4.26
$$ 
therefore we will identify from now on $\overline R$ to $\phi(\overline R)$. Since
$\overline f  \overline h = 2$, it follows that $\overline R$ is a conic bundle in
$\bold P^g$.  In particular $\overline R$ is a surface of degree $2g-4$ with
hyperelliptic hyperplane  sections. It will be called here \it Castelnuovo
surface \rm, (cfr. [C]). Let
$$
j: \overline R \to \overline R \tag 4.27
$$
be the involution induced by $\rho$, then $j$ acts on $H^0(\Cal O_{\overline
R}(\overline h))$.  The $+1$ eigenspace is $\rho^*H^0(\Cal O_R(h))$. The $-1$
eigenspace is generated by a section vanishing on $\rho^{-1}(E)$. Then
$\rho^*H^0(\Cal O_R(h))$ has codimension $1$ in $H^0(\Cal O_{\overline R}(\overline
h))$, moreover the map $\rho$ is induced by the projection from a point $o \in
\bold P^g$. One can also show that $\overline R$ is the intersection  of the cone
over $R$ of vertex $o$ with a quadric hypersurface. 
\bigskip \noindent (2) \it Pencils on the Castelnuovo surface. \rm
\par \noindent  Let $\Cal A \in Q_E$. For each $a \in \mid \Cal A \mid$, the linear
span $<2a>$ is a tangent hyperplane to the quadric cone $q$ such that $\epsilon_2(q)
= \Cal A$. This follows from the definition of~$\epsilon_2$.
\par \noindent (4.28) {\it CLAIM The curve $H_a = <2a> \cdot R$ is irreducible
for a general
$a$.} \par \noindent
{\it Proof}  Assume $H_a$ is reducible for each $a$. Then $H_a
=b_a + \Sigma f_{ia}$, where $b_a$ is an irreducible section of $R$ and $f_{ia} \sim
f$. Note that
$\mid a \mid$ is base-point-free because $E$ is general as in lemma (4.14). Then,
applying to $\mid a \mid$ the uniform position lemma, it turns out that  $a = \Sigma
x_i$, with $\lbrace x_i \rbrace = b_a \cap f_{ia}$. But then $f_{ia}\cdot E = 2x_i$
and $x_i$ is a Weierstrass point of $E$. This is impossible for a general $a$.
\bigskip \noindent  By the claim $H_a$ is a smooth rational curve. Since $H_a\cdot E =
2a$, we have~a~splitting 
$$
\rho^*H_a = A_a + j^*A_a, \tag 4.29
$$  
where $\rho/A_a: A_a \to H_a$ is biregular. By adjunction formula we have
$H_aK_R=-g$. Since $A_a\rho^{-1}(E)= H_aE = g-2$, it follows $A_aK_{\overline R}=-2$
and $A_a^2=0$. Moreover $A_a$ is not linearly equivalent to $j^*A_a$. Since
$R$ is regular, the curves of the family $ \Cal F = \lbrace A_a, j^*A_a / a \in \mid 
\Cal A \mid \rbrace$ belong to finitely many linear equivalence classes. 
One can show easily from this that $\Cal F$ is the disjoint union of two pencils:
$\mid A \mid$ and $\mid j^*A \mid$. From the exact sequence
$$
0 \to H^0(\Cal O_{\overline R}) \to H^0(\Cal O_{\overline R}(A)) \to H^0(\Cal
O_{A}(A))
\to 0
$$ 
it follows that $\mid A\mid$ is an irreducible base-point-free pencil on $\overline
R$.  We define as
$$
P_E  
$$ 
the set of all unordered pairs of line bundles $(\Cal O_{\overline R}(A), \Cal
O_{\overline R}(j^*A))$ which are constructed as above from some $\Cal A \in Q_E$. Let
$$
\epsilon_3: Q_E \to P_E \tag 4.30
$$
be the map sending $\Cal A$ to $(\Cal O_{\overline R}(A), \Cal O_{\overline
R}(j^*A))$.$\epsilon_3$ is surjective by definition of $P_E$. To show the
injectivity let $\Cal A, \Cal A'$ be distinct elements of $Q_E$ and let
$\epsilon_3(\Cal A)$ $=$ $(\Cal O_{\overline R}(A), \Cal O_{\overline R}(j^*A))$,
$\epsilon_3(\Cal A') = (\Cal O_{\overline R}(A'), \Cal O_{\overline R}(j^*A'))$.
The long exact sequence of
$$
0 \to \Cal O_{\overline R}(-E+A-A') \to \Cal O_{\overline R}(A-A') \to \Cal A \otimes
\Cal A^{'*} \to 0
$$
gives $h^0(\Cal O_{\overline R}(A-A')) = 0$, so that the two pairs are
distinct. $\epsilon_3$ is bijective. 
\bigskip \noindent (3) \it The surface $\tilde R$. \rm \par \noindent At first we
consider the curve $\overline C = \rho^*(C)$. This is a birational model of
$\tilde C$ with singular locus $\rho^*(d)$, where $2d = C \cdot E$. We can assume
that $Sing \overline C$ is a set of $2g-2$ ordinary nodes 
$$
\lbrace x_1, \dots , x_{2g-2} \rbrace 
$$
such that $2\Sigma \rho(x_i) = C \cdot E = 2d$. Let
$$
\sigma: \tilde R \to \overline R \tag 4.31
$$ be the blowing up along $\lbrace x_1, \dots, x_{2g-2} \rbrace$. Then the strict
transform of $\overline C$ is $\tilde C$, moreover the involution $j/\overline C:
\overline C \to \overline C$ lifts to the fixed-point-free involution $i$ of $\tilde
C$. Up to the action of $i^*$, each $(\Cal O_{\overline R}(A), \Cal
O_{\overline R}(j^*A)) \in P_E$ uniquely defines a pair
$$ (\tilde L_q, \Gamma_q), \tag 4.32
$$ 
where $\Gamma_q $ is the image of the restriction map $H^0(\Cal O_{\tilde
R}(\sigma^*A)) \to H^0(\Cal O_{\tilde C}(\sigma^*A))$ and~$\tilde L_q$~$\cong \Cal
O_{\tilde C}(\sigma^*A)$. Here the index $q$ denotes the quadric cone
corresponding to the element $(\Cal O_{\overline R}(A),\Cal O_{\overline R}(j^*A))$
under the bijection $\epsilon_3 \cdot
\epsilon_2$. \par \noindent
>From $\rho_* \sigma_*(\sigma^*A) = \rho_*A \in \mid h \mid$ it follows
$Nm \tilde L_q \cong \Cal O_C(\rho_*A) \cong \omega_C$. Let $\tilde l = \sigma^*A
\cdot \tilde C$, by projection formula one has $\rho_* \sigma_*(\sigma^*A \cdot
\tilde C)$ $=$ $\rho_*(A \cdot \sigma_*\tilde C) = \rho_*A \cdot C$ $=$
$\pi_* \tilde l$. On the other hand $\rho_*A \cdot E = 2a$ for some $a \in \mid \Cal A
\mid$. Since
$$
<\pi_* \tilde l> = <2a> = <\rho_*A>,
$$
the family $\lbrace <\pi_*\tilde l>, \tilde l \in \mid \Gamma_q \mid \rbrace$ is
exactly the family of tangent hyperplanes to $q$. Hence
$$
\alpha(\tilde L_q,\Gamma_q) = q.
$$
Let $K$ be a canonical divisor on $\tilde R$ and let
$$  u: \tilde R \to  \bold P^{2g-2} \tag 4.33
$$  be the map defined by $\mid K+\tilde C \mid$. Then $u/\tilde C$ is the canonical
embedding and we can put $\bold P = \bold P^{2g-2}$.  We consider a pair  $(\tilde
L_q, \Gamma_q)$ as above: for a divisor 
$\tilde l \in \mid \Gamma_q \mid$ we have
$$
\tilde l = \tilde C \cdot \sigma^*A \tag 4.34
$$  for some $A_{\tilde l}$ in the pencil $\mid A \mid$. 
\proclaim{4.35 LEMMA} $h^0(\tilde L_q) = h^0(\Cal O_{\tilde
R}(K+\tilde C-\sigma^*A))$.
\endproclaim
\demo {Proof} Recall that $h^1(\tilde L_q) = h^0(\tilde L_q)$ and
consider the standard exact~sequence
$$ 0 \to \Cal O_{\tilde R}(-\tilde C+\sigma^*A) \to \Cal O_{\tilde
R}(\sigma^*A) \to
\tilde L_q
\to 0.
$$ Since $h^0(-\tilde C+\sigma^*A) = h^1(\sigma^*A) = 0$, it follows
$h^1(\tilde L_q) = 2 + h^1(\Cal O_{\tilde R}(K+\tilde C- \sigma^*A))$. Now
consider the exact sequence
$$  0 \to \Cal O_{\tilde R}(K+\tilde C-\sigma^*A)) \to \Cal O_{\tilde R}(K+\tilde C)
\to \Cal O_{\sigma^*A}(K+\tilde C) \to 0.
$$  
We have $\Cal O_{\sigma^*A}(K+\tilde C) \cong \Cal O_{\bold P^1}(2g-4)$ and 
$\chi (\Cal O_{\tilde R}(K+\tilde C)) = 2g-1$, so that
$$
\chi(\Cal O_{\tilde R}(K+\tilde C-\sigma^*A)) = \chi (\Cal O_{\tilde R}(K+\tilde
C))-\chi (\Cal O_{\sigma^*A}(K+\tilde C)) = 2.
$$
This implies $h^0(\Cal O_{\tilde R}(K+\tilde C-\sigma^*A)) = 2 + h^1(\Cal
O_{\tilde R}(K+\tilde C - \sigma^*A))$ and the statement. \enddemo \noindent
Let $E' \subset \tilde R$ be the strict transform of $E$ by $\sigma$, we want to point
the equality 
$$
\mid K+\tilde C \mid = \mid E' + \sigma^*\overline h \mid. \tag 4.36
$$
Indeed it holds  $K \sim -\sigma^*\overline h+(g-4)\sigma^*\overline f + G$,
$G$ being the exceptional divisor of $\sigma$. Moreover $E' \sim \overline
\sigma^*h - G$ and $\tilde C \sim 3\sigma^*\overline h - (g-4)\overline \sigma^*f
- 2G$. In particular we have $K+\tilde C \sim 2\sigma^*\overline h - G \sim
E'+\sigma^*\overline h$. Notice also that 
$$\Cal O_{E'}(K+\tilde C) \cong \Cal O_E(2h-d) \tag 4.37,$$ with $2d = C\cdot E =
2\Sigma \rho(x_i)$. 
\proclaim {4.38 LEMMA} $h^0(\tilde L_q) = 2 + h^0(\Cal O_{E'}(K+\tilde
C-\sigma^*A))$. \endproclaim
\demo {Proof} The previous lemma implies $h^0(\tilde L_q)$ $=$ $h^0(\Cal O_{\tilde
R}(K+\tilde C-\sigma^*A))$. From (4.36) and
$A + j^*A \sim \overline h$ we have $K+\tilde
C-\sigma^*A-E'$ $\sim$ $\sigma^*j^*A$. Since $\mid j^*A \mid$ is a base-point-free
pencil of irreducible rational curves, the same holds for $\mid \sigma^*j^*A \mid$.
In particular we have $h^0(\Cal O_{\tilde R}(\sigma^*j^*A))=2$ and $h^1(\Cal
O_{\tilde R}(\sigma^*j^*A))=0$. Then the statement follows from the long exact 
sequence of
$$ 0 \to \Cal O_{\tilde R}(\sigma^*j^*A) \to \Cal O_{\tilde R}(K+\tilde
C-\sigma^*A)
\to \Cal O_{E'}(K+\tilde C-\sigma^*A) \to 0.
$$
\enddemo \noindent 
(4) \it Completion of the proof of theorem 4.19. \rm \par \noindent By (4.37) $\Cal
O_{E'}(\tilde C + K - \sigma^*A_a)$ $\cong$ $\Cal O_E(2h-d)\otimes \Cal A^*$. Since
$\Cal A^{\otimes 2}
\cong \Cal O_E(h)$,  it follows~$\Cal O_{E'}(\tilde C + K - \sigma^*A)$ $\cong \Cal
A^{\otimes 3}(-d) \cong \theta_q$, where $\theta_q$ is the theta characteristic
considered in (4.7) and $\epsilon_2(q) = \Cal A$. By lemma (4.38), $ h^0(\tilde L_q) =
2 + h^0(\theta_q)$.
\par \noindent
The motivation for studying the trigonal case was the proof of the next theorem.
\proclaim {4.39 THEOREM} Generically the function $b_+-b_-: \Cal R_g \to \bold
Z$ has value $2^{g-2}$. \endproclaim
\demo {Proof} By 3.21 it suffices to compute $b_+-b_-$ at a point of the set
$R$ defined in 3.10. By 4.9 an \'etale double covering of a general trigonal curve 
$C$ defines a point in $R$. We have seen that in this special case
$b_+(\pi)-b_-(\pi)$ is twice the difference between the number of even and 
odd thetas on a curve of genus $g-3$, that is $2^{g-2}$.
\enddemo \bigskip \noindent
\bf 5. Conclusion \rm \par \noindent
\proclaim {5.1 THEOREM} The degree of the Gauss map for the Theta divisor
of a general Prym variety $P$ is $D(g)+2^{g-3}$.
\endproclaim
\demo {Proof} Let $\pi$ be the \'etale double covering defining $P$.
Then $\pi$ is general and $b_+(\pi)$ is the degree of the Gauss map. 
By theorems 3.8 and 4.39 we have $b_+(\pi)+b_-(\pi) = 2D(g)$ and 
$ b_+(\pi)-b_-(\pi) = 2^{g-2}$. Hence $b_+(\pi) = D(g)+2^{g-3}$.
\enddemo \noindent
\bf 6. References \rm \par \noindent [ACGH]  E. Arbarello, M. Cornalba, P.A.
Griffiths, J. Harris \it Geometry of Algebraic Curves I, \rm  Springer New York
(1985) \par \noindent 
 [AM]  A. Andreotti, A. Mayer \it On period relations for abelian integrals on
algebraic curves, \rm  Ann. Scuola Norm. Sup. Pisa v.21 p.189-238 (1952)
\par \noindent 
[B] A. Bertram \it An existence theorem for Prym special divisors, \rm Invent. Math.
90 p.669-671 (1987) \par \noindent 
[B1]  A. Beauville \it Prym Varieties: A Survey \rm AMS Proceedings of
Symposia in Pure Mathematics,  v. 49 p.607-619 (1989)\par \noindent [B2] A. Beauville
\it Sous-vari\'et\'es speciales des vari\'et\'es de Prym, \rm  Compositio Math. v.45
p.357-383 (1982)
\par \noindent  [B3] A. Beauville \it Prym varieties and the Schottky problem,  \rm
Invent. math. v.41 p.149-196 (1977) \par \noindent
[C] G. Castelnuovo \it Sulle superfici algebriche le cui sezioni iperpiane sono curve
iperellittiche, \rm Rend. Circ. Matem. Palermo
 v.4 p.73-88 (1890) \par \noindent
[CD] F.Cossec, I.Dolgachev \it Enriques
Surfaces I, \rm Birkh\"auser Basel  (1989) \par \noindent 
[D] O. Debarre \it   Sur le probl\'eme de Torelli pour les vari\'et\'es de
Prym, \rm  Am. J. of Math. v. 111 p.111-134 (1989) \par \noindent [HT] J.Harris
L.W.Tu \it On symmetric and skew-symmetric determinantal varieties, \rm Topology v.23
p.71-84 (1984) \par
\noindent
 [K]  G. Kempf , \it On the geometry of a theorem of Riemann, \rm 
Annals of Math. v.98 p.178-185 (1973) \par \noindent 
[M] D. Mumford \it Prym varieties I, \rm Contributions to Analysis,
Academic Press  New York p.325-350 (1974)
\par \noindent [R] S.P. Recillas \it Jacobians of curves with $g^1_4$'s are the
Prym's of trigonal curves, \rm Bol. de la Soc. Mat. Mexicana v.19 p.9-13 (1974)
\par \noindent  [Ro] T.G. Room \it The Geometry of determinantal loci, \rm Cambridge
University Press Cambridge (1938)
\par \noindent  [T] A.N. Tjurin \it The geometry of the Poincar\'e theta-divisor of a
Prym variety, \rm  Math. USSR Izvestija v.9 p.951-986 (1975) \par \noindent [W] G.
Welters \it A theorem of Gieseker-Petri type for Prym varieties, \rm Ann. Sci. Ecole
Norm. Sup. v.18 p.671-683 (1985) \bigskip \noindent AUTHOR's ADDRESS: \par \noindent
Dipartimento di Matematica, Universita' di Roma Tre \par \noindent L.go S.Leonardo
Murialdo 1 \par \noindent 00146 ROME (Italy) \par \noindent e-mail:
verra\@matrm3.mat.uniroma3.it
\end